\documentclass[12pt]{amsart}
\usepackage{amsmath,xypic}
\usepackage{pinlabel}
\usepackage{amsfonts,amssymb}
\usepackage[all]{xy}
\theoremstyle{plain}
\newtheorem{theorem*}{Theorem}
\newtheorem*{lemma*}{Lemma}
\newtheorem{corollary*}{Corollary}
\newtheorem*{proposition*}{Proposition}
\newtheorem{conjecture*}{Conjecture}
\newtheorem{theorem}{Theorem}[section]

\newtheorem{conjecture}[theorem]{Conjecture}
\newtheorem{question}[theorem]{Question}

\theoremstyle{remark}

\newtheorem*{remark}{Remark}

\theoremstyle{definition}

\def\scrk{\mathcal{K}}
\def\scrn{\mathcal{N}}
\def\scra{\mathcal{A}}
\def\sl{\operatorname{SL}}

\def\gl{\mbox{GL}}  \def\F{\Bbb{F}} \def\Z{\Bbb{Z}} \def\R{\Bbb{R}} \def\C{\Bbb{C}}
\def\N{\Bbb{N}}    
 \def\a{\alpha}   \def\bp{\begin{pmatrix}}
\def\sm{\setminus} \def\ep{\end{pmatrix}} \def\bn{\begin{enumerate}} 
   \def\en{\end{enumerate}}
\def\ba{\begin{array}} \def\ea{\end{array}}  
   \def\a{\alpha}  \def\ti{\tilde}
\def\id{\mbox{id}}   
  
\def\be{\begin{equation}} \def\ee{\end{equation}}

     \def\fr12{\frac{1}{2}} \def\z12{\Z[\fr12]}

\def\tpm {[t^{\pm 1}]}

\def\G{\Gamma}

\def\vol{\operatorname{vol}}
\def\scrk{\mathcal{K}}

\def\G{\Gamma}
\def\ct{\C\tpm}

\def\T{\mathcal{T}}
\def\tkt{\T_K(t)}

\def\scrn{\mathcal{N}}

\begin{document}

\title{Approximations to the volume  of hyperbolic knots}

\author{Stefan Friedl}
\address{Mathematisches Institut\\ Universit\"at zu K\"oln\\   Germany}
\email{sfriedl@gmail.com}

\author{Nicholas Jackson}
\address{University of Warwick, Coventry, UK}
\email{nicholas.jackson@warwick.ac.uk}

\date{\today}
\begin{abstract}
We present computational data and heuristic arguments which suggest that given a hyperbolic knot the volume correlates with its  determinant, the Mahler measure of its Alexander polynomial and the Mahler measure of the twisted Alexander polynomial corresponding to the discrete and faithful $\sl(2,\C)$ representation.
\end{abstract}
\maketitle

\section{Introduction}

A 3-manifold is  called \emph{hyperbolic} if the interior of $N$ admits a complete metric of constant curvature $-1$ of finite volume.
If $N$ is hyperbolic, then it follows from Mostow-Prasad rigidity that the hyperbolic metric is unique up to isometry. In particular the volume of $N$ with respect to the hyperbolic metric is an invariant of the
3-manifold $N$, which we denote by $\vol(N)$. If $N$ is any irreducible, compact, orientable 3-manifold with empty or toroidal boundary,
then  we define $\vol(N)$ to be the sum of the volumes of the hyperbolic pieces in the JSJ decomposition of $N$.

We say that a knot  $K\subset S^3$ is hyperbolic, if $S^3\sm \nu K$ is hyperbolic. In that case we will denote the volume by $\vol(K)$.
Note that by Thurston's geometrization theorem any knot which is neither a torus knot nor a satellite knot, is a hyperbolic knot.

One of the most elementary invariants of a knot $K$ is its determinant $\det(K)$. The determinant has many equivalent definitions,
in particular it equals any of the following:
\bn
\item the evaluation at $t=-1$ of the Alexander polynomial $\Delta_K(t)$,
\item the evaluation at $t=-1$ of the Jones polynomial $J_K(t)$,
\item the order of the torsion part of the first homology of the 2-fold cyclic cover of $S^3\sm K$,
\item the order of the homology of the 2-fold branched cover of $K$.
\en
Computations by Dunfield \cite{Du99} showed that there exists a surprising correlation between the volume of a knot and its determinant. Other invariants of interest to us are the Alexander polynomial $\Delta_K(t)$, the Jones polynomial $J_K(t)$ and the symmetrized twisted Alexander polynomial $\tkt$ of $K$ corresponding to the discrete and faithful representation which was introduced in \cite{DFJ11}.
In this paper we will give computational and heuristic evidence  that the volume of a hyperbolic knot $K$ is  related to the Mahler measure of the Alexander polynomial, the evaluation at $t=-1$ and $t=1$ of $\tkt$ and
 the Mahler measure of $\tkt$.

\subsection*{Acknowledgment.} We would like to thank Nicholas Bergeron, Nathan Dunfield, Thang Le, Wolfgang L\"uck and Saul Schleimer  for helpful
conversations.

\section{Hyperbolic volume, $L^2$-invariants and the Alexander polynomial}\label{section:l2}

In this section we will give some heuristic arguments, using the theory of $L^2$-invariants, 
to justify why we think
that the determinant, the Mahler measure of the Alexander polynomial and the Mahler measure of $\tkt$ are related to the  volume of a hyperbolic knot $K$.

\subsection{$L^2$-invariants}\label{section:l2general}

We first recall a few properties of certain $L^2$-invariants.
We refer to L\"uck's monograph \cite{Lu02} for full details.
Let $f$ be one of the following `classical invariants':
\bn
\item the $i$-th Betti number of a topological space $X$,
\item the signature of a compact even dimensional manifold $X$,
\item the Atiyah-Patodi-Singer $\eta$-invariant of an odd-dimensional closed Riemannian manifold $X$  (see \cite{APS75} for details).
\en
Each of the above invariants admits an `$L^2$-version'.
More precisely, given $f$ and $X$ as above and a homomorphism $\varphi\colon\pi=\pi_1(X)\to G$ to a (not necessarily finite) group $G$ there exists a real valued invariant $f^{(2)}(X,\varphi)$ which has in particular the following three properties:
\bn
\item[(A)] Let $\varphi\colon\pi\to G$ be an epimorphism onto a finite group $G$. If we denote by $X_\varphi$  the cover of $X$ corresponding to $\varphi$, then
\[ f^{(2)}(X,\varphi)=\frac{1}{|G|}f(X_\varphi).\]
\item[(B)] Let $\G_n$ be a nested sequence of normal subgroups of $\pi$, i.e. a sequence of normal subgroups of the form
\[ \pi \supset \G_1\supset \G_2 \supset \G_3\dots \]
then
\[ f^{(2)}(X,\pi\to \pi/\cap \G_n)=\lim_{n\to \infty} f^{(2)}(X,\pi\to \pi/\G_n).\]
\item[(C)] If $\psi\colon G \to H$ is a monomorphism, then
\[ f^{(2)}(X,\varphi)=f^{(2)}(X,\psi\circ \varphi).\]
\en

\subsection{The $L^2$-torsion}

Let $C_*$ be a chain complex over a field $\F$, together with a basis $c_*$ for $C_*$.
If $C_*$ is acyclic, then its torsion $\tau(C_*,c_*)\in \F^\times=\F\sm\{0\}$ is defined.
We refer to \cite{Mi66} and \cite{Tu01} for the definition.
(Note that the convention used by Turaev \cite{Tu01} gives the multiplicative inverse of the torsion invariant defined by Milnor, we will follow Milnor's convention throughout.)
Now let $X$ be  a  finite CW-complex. The cells naturally give rise to a basis for the chain complex $C_*(X;\R)$.
But $H_0(X;\R)$ is always non-zero, which means that the based complex $C_*(X;\R)$ is not acylic, i.e. the torsion of $C_*(X;\R)$ is not defined.

This problem can be circumvented by appealing to a more elaborate version of torsion.
Indeed, let $C_*$ be a chain complex over a field $\F$, together with a basis $c_*$ for $C_*$ and a basis $h_*$ for the homology groups,  then a torsion-invariant  $\tau(C_*,c_*,h_*)\in \F^\times$ is defined.
(We again refer to \cite{Mi66} and \cite{Tu01} for details.)

Now let $X$ again be a finite CW-complex. The cells give again rise to a basis $c_*$ of $C_*(X;\R)$.
We pick a basis $h_*$ for $H_*(X;\Z)/\mbox{torsion}$. This basis gives rise to a basis for
$H_*(X;\R)$ which we also denote by $h_*$. A straightforward argument shows that
up to a sign the invariant $\tau(C_*(X;\R),c_*,h_*)$ is independent of the choice of $h_*$.
Note that torsion is a multiplicative invariant, to keep the analogy with the previous section we now define
\[ \tau(X):=\ln|\tau(C_*(X;\R),c_*,h_*)|,\]
which is an additive invariant.
By \cite[Theorem~4.7]{Tu01} this invariant can be computed as follows:
\[ \tau(X)=\sum_{i} (-1)^i \ln  |\operatorname{Tor}(H_i(X;\Z))|.\]
Note that if $X$ is a 3-manifold, then it follows from Poincar\'e duality that
\[ \tau(X)=-\ln  |\operatorname{Tor}(H_1(X;\Z))|.\]

The  awkwardness in the definition of the classical torsion of a finite CW-complex translates into serious technical difficulties for the $L^2$-torsion. In particular, given a finite CW-complex $X$ and a group homomorphism $\varphi$ the $L^2$-torsion
$\tau^{(2)}(X,\varphi)$ is in general \emph{not} defined, even if all $L^2$-Betti numbers vanish
(see \cite[Section~3]{Lu02} for details). A sufficient (but not necessary) condition for the $L^2$-torsion to be defined is that the $L^2$-Betti numbers vanish and the Novikov-Shubin invariants are positive.
(Note that there are two conventions for the $L^2$-torsion: The torsion as defined in \cite{Lu02} equals the torsion of \cite{LS99} \emph{divided by two}, we will follow the convention of \cite{Lu02}.)

A general approximation  result for $L^2$-torsion  has not been proved yet, in particular the following two questions (which are a variation on \cite[Question~13.73]{Lu02}) seem to be wide open:

\begin{question}\label{question:approxl2}
Let $X$ be a finite CW-complex and let   $\G_n$ be a nested sequence of normal subgroups of $\pi=\pi_1(X)$.
Suppose $\tau^{(2)}(X,\pi\to \pi/\cap \G_n)$ is defined and suppose that $\tau^{(2)}(X,\pi\to \pi/\G_n)$ is defined for any $n$,
does the following equality hold:
\[ \tau^{(2)}(X,\pi\to \pi/\cap \G_n)=\lim_{n\to \infty} \tau^{(2)}(X,\pi\to \pi/\G_n)\,\,\,?\]
\end{question}

\begin{question}\label{question:approxl2finite}
Let $X$ be a finite CW-complex and let   $\G_n$ be a nested sequence of finite index normal subgroups of $\pi=\pi_1(X)$.
We denote by $X_n$ the cover corresponding to $\G_n$.
Suppose that $\tau^{(2)}(X,\pi\to \pi/\cap \G_n)$,
does the following equality hold:
\[  \tau^{(2)}(X,\pi\to \pi/\cap \G_n)=\lim_{n\to \infty} \frac{1}{|X:X_n|} \tau(X_n)\,\,\,?\]
\end{question}

We refer to \cite{BV10} for a detailed discussion of many related questions. We also refer to \cite[Section~4]{Se10} for a helpful outline of the philosophy relating torsion and volume.

\subsection{$L^2$-torsion of knots and the Mahler measure}

Let  $K\subset S^3$ be a knot.
Throughout this section we denote by $X=S^3\sm \nu K$ the exterior of $K$ and we write $\pi=\pi_1(X)$. We will equip $X$ with the structure of a finite CW complex.
The only $L^2$-torsions which are well understood are the  $L^2$-torsions corresponding to the identity map $\id\colon\pi\to \pi$ and corresponding to the abelianization map $\pi\to \Z$ which we denote by $\a$.
More precisely the following holds:

\begin{theorem}\label{thm:lueck}
Let $K\subset S^3$ be a knot.  Then $\tau^{(2)}(X,\id)$ and
$\tau^{(2)}(X,\a)$ are defined and they are given as follows:
\[ \ba{rcl}   \tau^{(2)}(X,\id)&=&-\frac{1}{6\pi} \operatorname{vol}(K),\\[1mm]
\tau^{(2)}(X,\a)&=&-\ln(m(\Delta_K(t))),\ea \]
where $m(\Delta_K(t))$ denotes the Mahler measure of the Alexander polynomial (we refer to the appendix for the definition). 
\end{theorem}

We refer to \cite[p.~206~and~Equation~(3.23)]{Lu02} (see also \cite[Lemma~13.53]{Lu02})  for a proof regarding $\tau^{(2)}(S^3\sm K,\a)$
 and we refer to \cite[Theorem~4.3]{Lu02} for details regarding  $\tau^{(2)}(S^3\sm K,\id)$ (see also \cite[Theorem~0.7]{LS99}).
The second statement was also proved by Li and Zhang (see \cite[Equation~8.2]{LZ06}). Note that the first statement in fact holds for any irreducible 3-manifold.

The Mahler measure of the Alexander polynomial has been studied extensively.
In particular Silver and Williams \cite[Theorem~2.1]{SW02} proved the following theorem
which gives an affirmative answer to Question \ref{question:approxl2finite} in a special case.

\begin{theorem} \label{thm:sw}
Let $K\subset S^3$ be a knot. Given $n\in \N$ we denote by $X_n$ the $n$-fold cyclic cover of $X=S^3\sm K$.
Then the following equality holds:
\[ \ln(m(\Delta_K(t)))=\lim_{n\to \infty} \frac{1}{n} \ln |\operatorname{Tor} H_1(X_n;\Z)|,\]
or equivalently
\[ \tau^{(2)}(X,\a)=\lim_{n\to \infty}\frac{1}{n}\tau(X_n).\]
\end{theorem}

\begin{remark}
\bn
\item
This theorem was proved by \cite{Ri90} and Gonz\'alez-Acu\~na and Short \cite{GS91} for the subsequence of $X_n$'s
for which the homology of the $n$-th fold branched cover is finite.
\item The theorem was recently reproved by Bergeron and Venkatesh \cite[Theorem~7.3]{BV10} and it was also recently extended by Raimbault \cite{Ra10} and Le \cite{Le10} to more general cases.
\item
Let $Y$ be any 3-manifold with empty or toroidal boundary. We denote by $\{Y_n\}$ the directed system of all finite covers of $Y$.
Thang Le \cite{Le09} has recently shown that  the following inequality holds:
\[ \lim \sup_{n\to \infty} \frac{1}{[Y:Y_n]}\ln |\operatorname{Tor} H_1(Y_n;\Z)|\leq \frac{1}{6\pi}\vol(Y),\]
or equivalently
\[ \lim \sup_{n\to \infty} \frac{1}{[Y:Y_n]}\tau(Y_n)\geq \tau^{(2)}(Y,\id).\]
\en
\end{remark}

\section{The hyperbolic torsion}

Given any orientable hyperbolic 3-manifold $Y$  we can consider the
discrete and faithful $\sl(2,\C)$--representation $\a_{can}$.
The corresponding twisted chain complex
$$C_*(\ti{Y})\otimes_{\Z[\pi_1(Y)]}\C^2$$
is acylic (see  \cite{Po97} and \cite{MP10}), and it follows that the
corresponding torsion $\tau(Y,\a_{can})\in \R\sm \{0\}$ is defined.

We recall the following  recent result of  Bergeron and Venkatesh (see \cite[Theorem~4.5]{BV10} and Example (3) of  \cite[Section~5.9.3]{BV10} with $(p,q)=(1,0)$).

 \begin{theorem}
 Let $Y$ be a closed hyperbolic 3-manifold. Let $\{Y_n\}_{n\in \N}$ be a nested collection of finite covers such that
\[ \bigcap\limits_{n\in \N} \pi_1(Y_n)\subset \pi_1(Y)\]
is trivial. Then the following holds:
\[ \lim_{n\to \infty} \frac{1}{[Y:Y_n]}\ln |\tau(Y_n,\a_{can}))|=-\frac{13}{6\pi} \vol(Y).\]
\end{theorem}

The question naturally arises, whether the conclusion of the theorem also holds for hyperbolic 3-manifolds with toroidal boundary.

 \begin{conjecture}
 Let $Y$ be a hyperbolic 3-manifold of finite volume. Let $\{Y_n\}_{n\in \N}$ be a nested collection of finite covers such that
\[ \bigcap\limits_{n\in \N} \pi_1(Y_n)\subset \pi_1(Y)\]
is trivial. Then there exists a constant $C>0$ such that
\[ \lim_{n\to \infty} \frac{1}{[Y:Y_n]}\ln |\tau(Y_n,\a_{can})|= -C\cdot \vol(Y)\]
for some $C$ independent of $Y$.
\end{conjecture}

Unfortunately at the moment there does not seem to be a good candidate for the constant $C$ in the conjecture.

Given a hyperbolic knot $K$ the authors and Nathan Dunfield introduced in  \cite{DFJ11}  an invariant $\T_K(t)\in \ct$ which  is defined as the normalized twisted Alexander polynomial of a hyperbolic knot corresponding to the discrete and faithful $\sl(2,\C)$ representation of the knot group.
It follows from the definition and standard arguments  (see e.g. \cite{DFJ11}) that the following equality holds:
Let $K$ be a hyperbolic knot and denote by $Y_n$ the $n$-fold cyclic cover of $Y=S^3\sm K$, then
\[ \lim_{n\to \infty} \frac{1}{n}\ln |\tau(Y_n,\a_{can}))|=-\ln(m(\T_K(t))).\]
In light of the above results it is therefore perhaps not surprising that the natural logarithm of the Mahler measure of $m(\T_K(t))$ correlates strongly with  $ \vol(Y)$.
The calculations of the next section below suggest that for knot complements the constant $C$ of the conjecture should be
\[ \lim_{n\to \infty} \frac{1}{[Y:Y_n]}\ln |\tau(Y_n,\a_{can}))|\approx \ln(m(\T_K(t)))\approx  -0.296509767136\cdot \vol(K).\]
Hence $C\approx  0.296509767136$, but the authors do not know a `natural' candidate for what $C$ should be.

\begin{remark}
Note that for any $m$ there exists a unique irreducible representation $\varphi_m:\sl(2,\C)\to \gl(m+1,\C)$ given by the action of $\sl(2,\C)$ on the symmetric powers of $\C^2$.
Let $Y$ be a closed hyperbolic 3-manifold and $\rho:\pi_1(Y)\to \sl(2,\C)$ a discrete and faithful representation.
M\"uller \cite{Mu09,Mu10}  showed that
\[  \vol(Y) =\lim_{m\to \infty} -\frac{1}{m^2} 4\pi \ln \tau(Y,\varphi_m\circ \rho).\]
It is unfortunately not clear to the authors how this beautiful result relates to the above open questions.
\end{remark}

\section{Calculations}

In this section we will compute various real valued invariants for all knots up to 15 crossings
and we will compare the values to the volume.

\subsection{Comparison of invariants}

Before we present our calculations we first introduce  some notation and definitions.
Given a set $\scrk$ of hyperbolic knots and a real-valued invariant $\Phi$ of hyperbolic knots
we denote by
\[ \ba{rcl} A(\Phi,\scrk)&:=&\frac{1}{|\scrk|} \sum_{K\in \scrk} \Phi(K),\\
\sigma(\Phi,\scrk)&:=&\sqrt{\frac{1}{|\scrk|} \sum_{K\in \scrk} \left(\Phi(K)-A(\Phi,\scrk)\right)^2}\ea \]
the average value respectively  the standard deviation of $\Phi$. When $\scrk$ is understood, then we will drop it from the notation.
We furthermore define
\[ \ba{rcl}  A_{vol}(\Phi,\scrk)&:=& A(\Phi/\vol,\scrk),\\
 \sigma_{vol}(\Phi,\scrk)&:=& \sigma(\Phi/\vol,\scrk),\\
  \Sigma_{vol}(\Phi,\scrk)&:=&\frac{1}{A_{vol}(\Phi,\scrk)}\cdot \sigma_{vol}(\Phi,\scrk).\ea \]
Finally we define $r(\Phi,\scrk)$ to be the Pearson correlation coefficient of the invariant $\Phi$ and hyperbolic volume, i.e.
\[ r(\Phi,\scrk)= \frac{1}{|\scrk|-1}\sum_{K\in \scrk}\frac{\Phi(K)-A(\Phi)}{\sigma(\Phi)}\cdot \frac{\vol(K)-A(\vol)}{\sigma(\vol)}.\]
Note that $r(\Phi,\scrk)=1$ if there exists an $s>0$ such that $\Phi(K)=s\vol(K)$ for all $K\in \scrk$,
and $r(\Phi,\scrk)=-1$ if there exists an $s<0$ such that $\Phi(K)=s\vol(K)$ for all $K\in \scrk$. In general $r(\Phi,\scrk)\in [-1,1]$. The absolute value of $r(\Phi,\scrk)$ can be seen as a measure of linear dependence between volume and $\Phi$. For all of the above we will drop the invariant $\Phi$ from the notation, if $\Phi$ is understood.

\subsection{The data from polynomials}

In the following pages we will plot the volume of a hyperbolic knot against various invariants.
The diagrams are drawn using a randomly chosen sample of one quarter of all hyperbolic knots with at most fifteen crossings.
In the diagrams the data for alternating knots are plotted in red, while those
for non-alternating knots are shown in green.

Given $n$ we denote by $\scra_n$ the set of all hyperbolic alternating knots up to $n$  crossings, we denote by $\scrn_n$ the set of all hyperbolic non-alternating knots with up to $n$ crossings
and we denote by $\scrk_n$ the set of all hyperbolic knots with up to $n$ crossings.

\begin{figure}[htbp!]
\centering
\labellist\footnotesize
\pinlabel {$\ln|\Delta_K(-1)|$} [b] at 240 45
\pinlabel {\rotatebox{90}{$\vol(S^3\setminus K)$}} [l] at 55 185
\endlabellist
\includegraphics[width=0.8\textwidth]{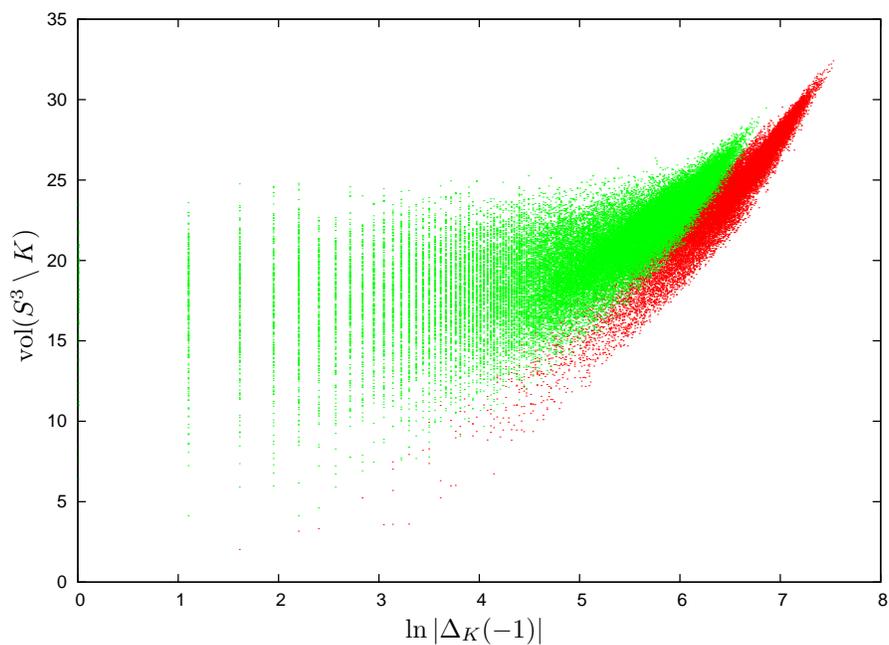}
\caption{Plot of $\vol(S^3\setminus K)$ against $\ln|\Delta_K(-1)|$
for hyperbolic knots $K$ with at most fifteen crossings.}
\label{fig:detvkc}
\end{figure}

\begin{figure}[htbp!]
\centering
\labellist\footnotesize
\pinlabel {$\ln(m(\Delta_K))$} [b] at 240 45
\pinlabel {\rotatebox{90}{$\vol(S^3\setminus K)$}} [l] at 55 185
\endlabellist
\includegraphics[width=0.8\textwidth]{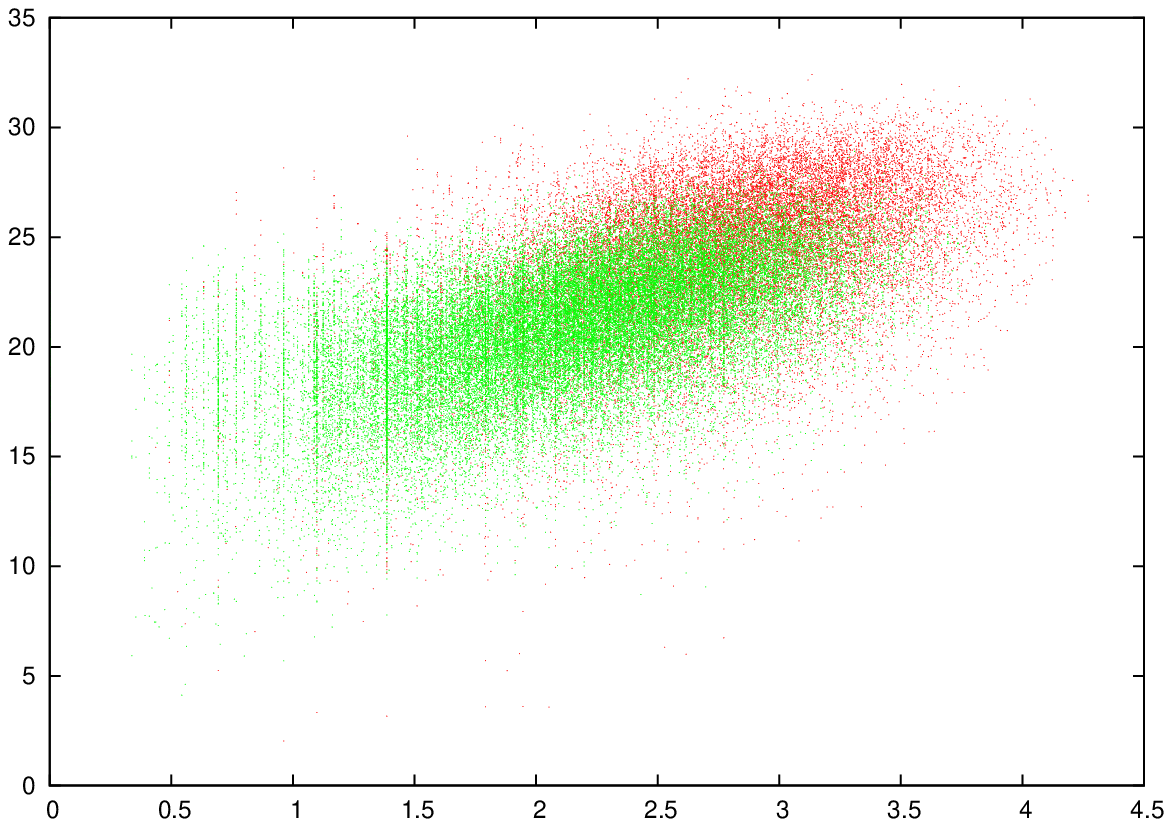}
\caption{Plot of $\vol(S^3\setminus K)$ against $\ln(m(\Delta_K))$
for hyperbolic knots $K$ with at most fifteen crossings.}
\label{fig:mmavkc}
\end{figure}

\begin{figure}[htbp!]
\centering
\labellist\footnotesize
\pinlabel {$\ln|\T_K(-1)|$} [b] at 240 45
\pinlabel {\rotatebox{90}{$\vol(S^3\setminus K)$}} [l] at 55 185
\endlabellist
\includegraphics[width=0.8\textwidth]{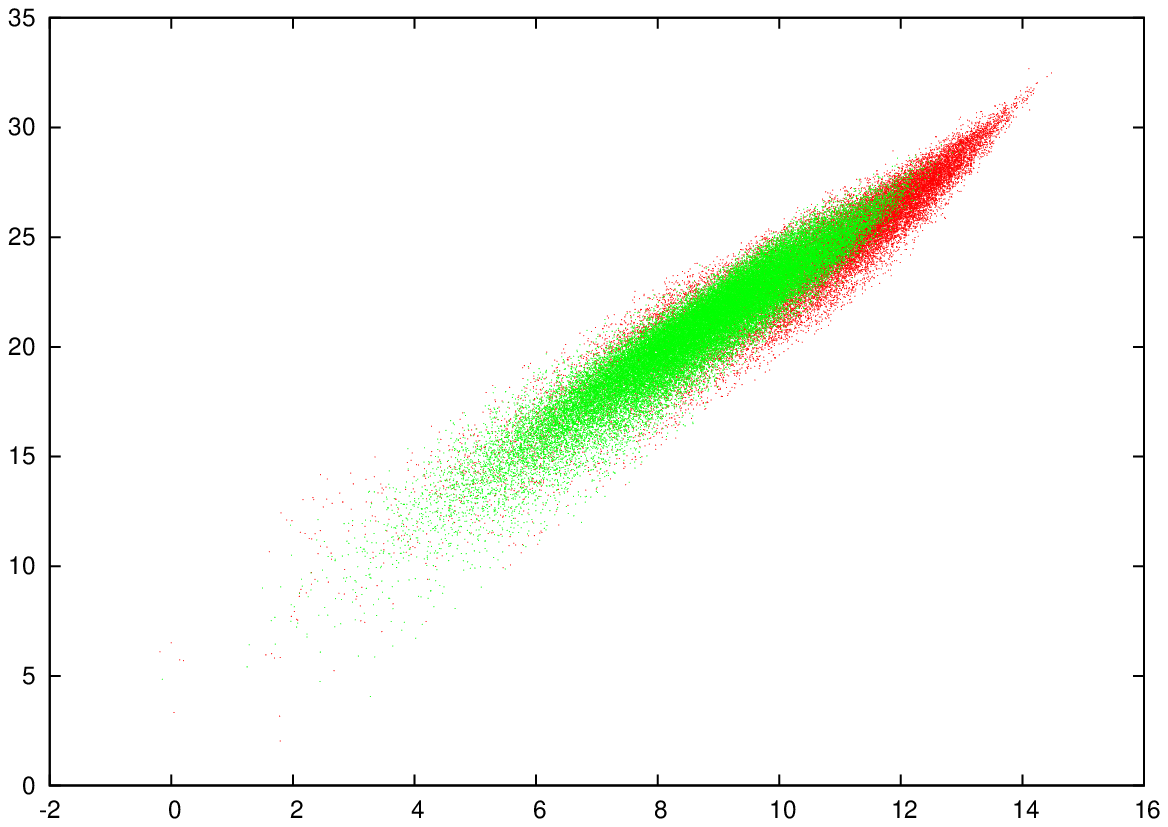}
\caption{Plot of $\vol(S^3\setminus K)$ against $\ln|\T_K(-1)|$ for
hyperbolic knots $K$ with at most fifteen crossings.}
\label{fig:temvkc}
\end{figure}

\begin{figure}[htbp!]
\centering
\labellist\footnotesize
\pinlabel {$\ln|\T_K(+1)|$} [b] at 240 45
\pinlabel {\rotatebox{90}{$\vol(S^3\setminus K)$}} [l] at 55 185
\endlabellist
\includegraphics[width=0.8\textwidth]{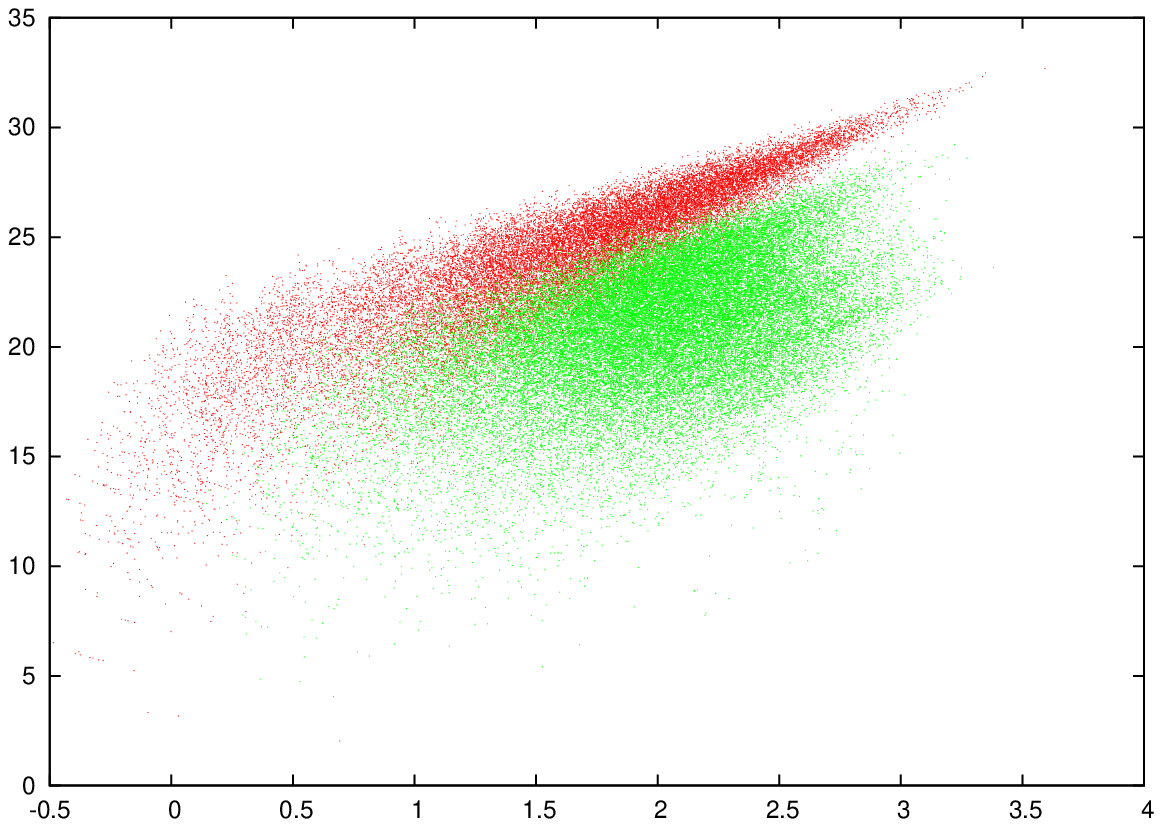}
\caption{Plot of $\vol(S^3\setminus K)$ against $\ln|\T_K(+1)|$ for
hyperbolic knots $K$ with at most fifteen crossings.}
\label{fig:tepvkc}
\end{figure}

\begin{figure}[htbp!]
\centering
\labellist\footnotesize
\pinlabel {$\ln(m(\T_K))$} [b] at 240 45
\pinlabel {\rotatebox{90}{$\vol(S^3\setminus K)$}} [l] at 55 185
\endlabellist
\includegraphics[width=0.8\textwidth]{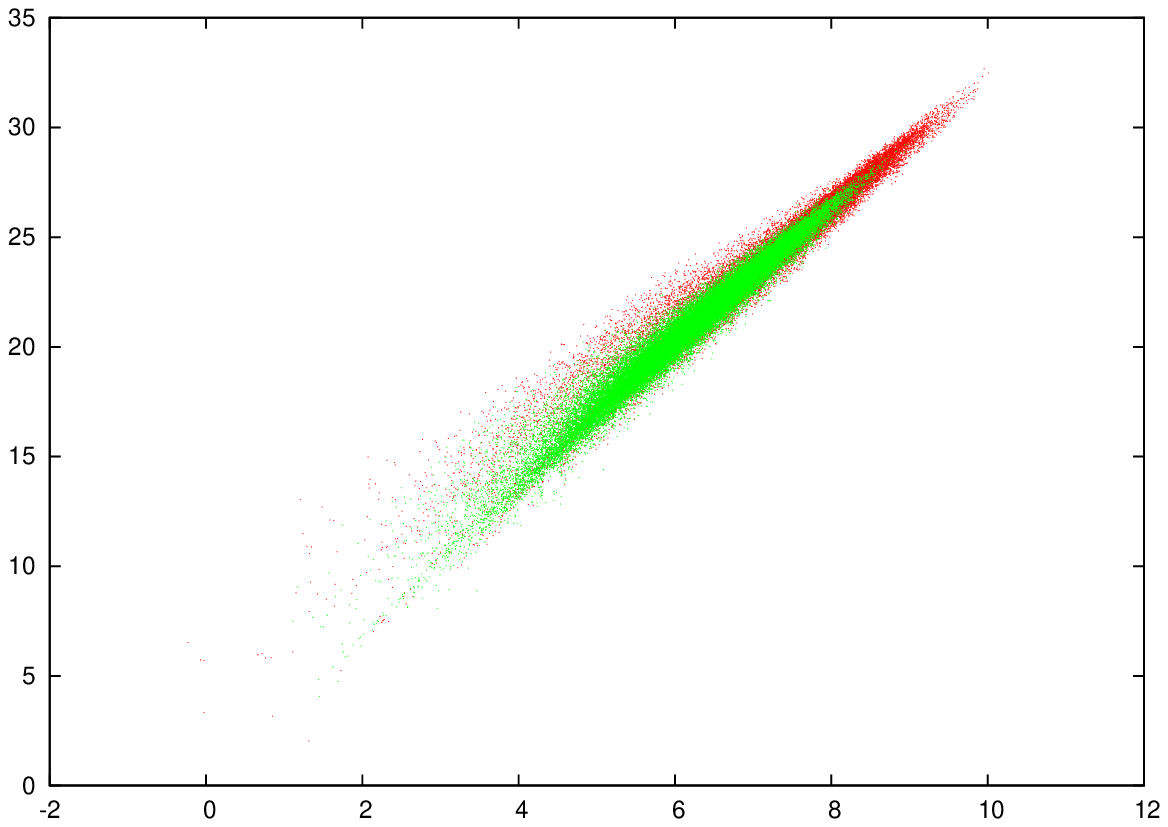}
\caption{Plot of $\vol(S^3\setminus K)$ against $\ln(m(\T_K))$ for
hyperbolic knots $K$ with at most fifteen crossings.}
\label{fig:mmtvkc}
\end{figure}

\begin{figure}[htbp!]
\centering
\labellist\footnotesize
\pinlabel {$\ln(m(J_K))$} [b] at 240 45
\pinlabel {\rotatebox{90}{$\vol(S^3\setminus K)$}} [l] at 55 185
\endlabellist
\includegraphics[width=0.8\textwidth]{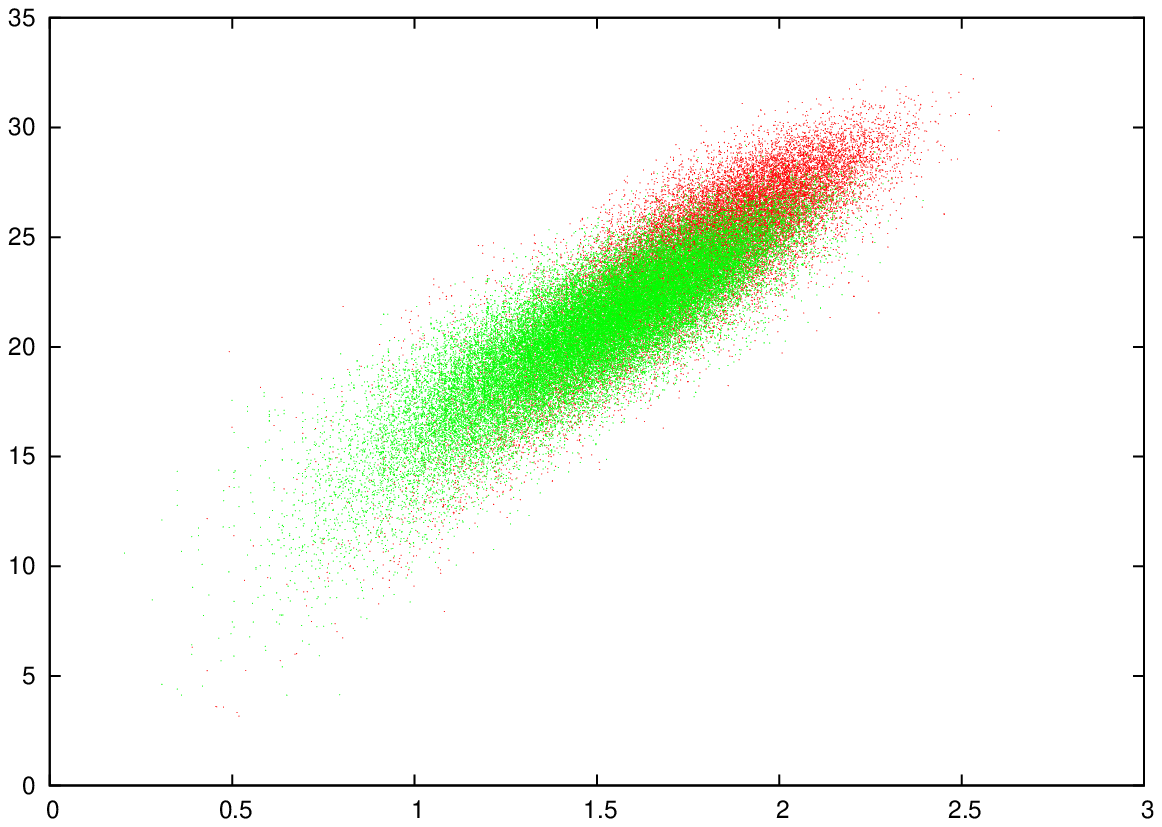}
\caption{Plot of $\vol(S^3\setminus K)$ against $\ln(m(J_K))$ for
hyperbolic knots $K$ with at most fifteen crossings.}
\label{fig:mmjvkc}
\end{figure}

The calculations were done with $n=15$, and the results are shown in
Tables~\ref{tbl:data1}--\ref{tbl:data3}.  Note that all but 22 prime
knots with 15 crossings or less are hyperbolic (see \cite{HTW98}).

\begin{table}[htbp!]
\centering
\footnotesize
\begin{tabular}{|r|l|l|l|l|l|} \hline
Invariant & $A_{vol}(\scra_{15})$ &  $\sigma_{vol}(\scra_{15})$ & $\Sigma_{vol}(\scra_{15})$
& $r(\scra_{15})$ \\
\hline
 $ \ln|\Delta_K(-1)|$ &
 0.271768541 & 0.0244239958 & 0.0898705780 &  0.9636372551 \\[1mm]
$\ln(m(\Delta_K(t)))$ &
 0.105675594& 0.0234634945& 0.2062259213 &  0.4779954145 \\[1mm]
\hline
$\ln(|\T_K(-1)|)$ &
  0.4400177279 & 0.0322204235 & 0.0732252850 &  0.9426517821 \\[1mm]
$\ln(|\T_K(+1)|)$ & 0.0643061411 & 0.0213629175 & 0.3322064909 &
 0.8930977985 \\
$\ln(m(\T_K(t)))$ &
  0.2975019481 & 0.0152639207 & 0.0513069604  & 0.9782289355 \\[1mm]
\hline
$\ln(|m(J(t))|)$ &
  0.0735034859 & 0.0063649487 & 0.0865938353 & 0.8465086985 \\[1mm] \hline
\end{tabular}
\medskip
\caption{Calculated data for alternating knots with up to fifteen
crossings}
\label{tbl:data1}
\end{table}

\begin{table}[htbp!]
\footnotesize
\begin{tabular}{|r|l|l|l|l|l|} \hline
Invariant & $A_{vol}(\scrn_{15})$ &  $\sigma_{vol}(\scrn_{15})$ & $\Sigma_{vol}(\scrn_{15})$
  &  $r(\scrn_{15})$ \\ \hline
$\ln|\Delta_K(-1)|$ &
 0.245020662 & 0.0947595379 & 0.3867410076 &  0.7048463196 \\
$\ln(m(\Delta_K(t)))$ &
 0.101165642 & 0.0259754182& 0.2567612649 &  0.5658607451 \\
  \hline
$\ln(|\T_K(-1)|)$ &
  0.4172095065 & 0.0276945418 & 0.0663804189 &  0.9475455803 \\
$\ln(|\T_K(+1)|)$ &
  0.0958788600 & 0.0211610565 & 0.2207061752 & 0.5216657082 \\
$\ln(m(\T_K(t)))$ &
  0.2959611229 & 0.0104234142 & 0.0352188629 & 0.9825263792 \\
  \hline
$\ln(|m(J(t))|)$ &
  0.0718864873 & 0.0074772898 & 0.1040152344 & 0.8422774659 \\ \hline
\end{tabular}
\medskip
\caption{Calculated data for non-alternating knots with up to fifteen
crossings}
\label{tbl:data2}
\end{table}

\begin{table}[htbp!]
\footnotesize
\begin{tabular}{|r|l|l|l|l|l|}
\hline
Invariant & $A_{vol}(\scrk_{15})$ &  $\sigma_{vol}(\scrk_{15})$ & $\Sigma_{vol}(\scrk_{15})$
&  $r(\scrk_{15})$ \\ \hline
$ \ln|\Delta_K(-1)|$ &
  0.25458698& 0.0783940759 & 0.3079264927 &  0.7782203193 \\
$\ln(m(\Delta_K(t)))$ &
 0.105675594 & 0.0258232111 & 0.2443630571 & 0.6210850376 \\
\hline
$\ln(|\T_K(-1)|)$ &
  0.4253305875 & 0.0313498808 & 0.0737070922 & 0.9558964834 \\
$\ln(|\T_K(+1)|)$ &
  0.0846370975 & 0.0260653503 & 0.3079660223 &  0.4394308257 \\
$\ln(m(\T_K(t)))$ &
  0.2965097482 & 0.0123880680 & 0.0417796315 &  0.9834187657 \\
  \hline
$\ln(|m(J(t))|)$ &
  0.0724622238 & 0.0071433263 & 0.0985800043 & 0.8711042342 \\ \hline
\end{tabular}
\medskip
\caption{Calculated data for all knots with up to fifteen crossings}
\label{tbl:data3}
\end{table}

%
%

We would like to point out several observations:
\bn
\item Note that the average ratio of $m(\Delta_K(t))$ to volume is
\[  0.105675594\cong  \frac{1}{3\pi} \cdot 0.995969009.\]
At least a na\"{\i}ve reading of Section \ref{section:l2} would have suggested that the ratio should on average be $\frac{1}{6\pi}$.
\item According to the Pearson correlation coefficient the determinant of $K$ correlates better with the hyperbolic volume than the Mahler measure of the Alexander polynomial, especially for alternating knots.
This also seems somewhat surprising considering the discussion of Section \ref{section:l2}.
\item The Mahler measure of $\tkt$ correlates very highly with the volume. The average ratio of $0.2965097482 $ is rather mysterious though, it is not clear what `nice' number it corresponds to.
\en

We conclude with a list of open questions:
\bn
\item Does the Mahler measure of knot Floer homology (viewed as a multivariable polynomial) correlate more with  hyperbolic volume than the Mahler measure of the Alexander polynomial?
\item Does the Mahler measure of Khovanov homology (viewed as a multivariable polynomial) correlate more with hyperbolic volume than the Mahler measure of the Jones polynomial?
\item Is there any correlation between the Chern-Simons invariant and the arguments of the zeros of the Alexander polynomial?
\en

\section{Appendix: The Mahler measure}

In this appendix we will quickly recall the definition and basic properties of the Mahler measure of a polynomial.
Throughout this section we refer to \cite[Section~16]{Sc95} and \cite{SW04} for details and references.

Let $p\in \C[t_1^{\pm 1},\dots,t_m^{\pm 1}]$ be a multivariable Laurent polynomial.
The \emph{Mahler measure $m(p)$} of a non-zero polynomial is defined to be
\[ m(p)=\exp \int_{(S^1)^m} \ln|p(s)|ds \]
where we equip the $m$-torus $(S^1)^d$ with the Haar measure.
(Note that the integral is well-defined despite the zeros of $p$.)

In the case of one-variable polynomials we can rewrite the definition as follows:
Let $p(t)\in \ct$ be a polynomial. Then
\[ m(p(t))=\exp \frac{1}{2\pi} \int_{\theta=0}^{2\pi} \ln |p(e^{i\theta})| d\theta.\]
If
\[ p(t)=ct^k\prod_{i=1}^n (t-r_i),\]
then it follows from  Jensen's formula that
\[ m(p(t)) = |c| \cdot \prod\limits_{j=1}^n \operatorname{max}(|r_j|,1).\]


\begin{thebibliography}{asfd}

\bibitem[APS75]{APS75}
M. Atiyah, V. Patodi and  I. Singer, {\em Spectral asymmetry and Riemannian geometry II}, Math. Proc. Camb. Phil. Soc, 78:
405-432
(1975)

\bibitem[BV10]{BV10}
N. Bergeron and A. Venkatesh, {\em The asymptotic growth of torsion homology for arithmetic groups}, Preprint (2010)

\bibitem[Du99]{Du99}
N. Dunfield, {\em An interesting relationship between the Jones polynomial and hyperbolic volume}, unpublished note (1999)\\
\texttt{http://www.math.uiuc.edu/\,\~\,nmd/preprints/misc/dylan/index.html}

\bibitem[DFJ11]{DFJ11}
N. Dunfield, S. Friedl and N. Jackson, {\em Twisted Alexander polynomials of hyperbolic knots},
in preparation (2011)

\bibitem[GS91]{GS91}
F. Gonz\'alez-Acu\~na and H. Short, {\em Cyclic branched coverings of knots and homology
spheres}, Revista Math. 4 (1991), 97 - 120.


\bibitem[HTW98]{HTW98}
J. Hoste, M. Thistlethwaite and J. Weeks, {\em The first 1,701,936 knots}, Math. Intell. 20(4) (1998), 33--48.




\bibitem[Le09]{Le09}
T. Le, {\em Hyperbolic volume, Mahler measure, and homology
growth}, talk at Columbia University (2009), slides available from\\
\texttt{http://www.math.columbia.edu/\,\~\,volconf09/notes/leconf.pdf}

\bibitem[Le10]{Le10}
T. Le, {\em
Homology torsion growth and Mahler measure}, Preprint (2010)


\bibitem[LZ06]{LZ06}
W. Li and W. Zhang,  {\em An $L^2$-Alexander invariant for knots}, Commun. Contemp.
Math., 8(2):167–187, 2006.



\bibitem[L\"u94]{Lu94} W. L\"uck, {\em
Approximating $L^2$-invariants by their finite-dimensional analogues},
Geom. Funct. Anal., 4(4):455-481, 1994.

\bibitem[LS99]{LS99}
W. L\"uck and  T. Schick, {\em $L\sp 2$-torsion of hyperbolic manifolds of finite volume}, Geom.
Funct. Anal.  9  (1999),  no. 3, 518--567.

\bibitem[L\"u02]{Lu02} W. L\"uck, {\em
 $L\sp 2$-invariants: theory and applications to geometry and $K$-theory}, Ergebnisse
der Mathematik und ihrer Grenzgebiete. 3. Folge. A Series of Modern Surveys in Mathematics, 44.
Springer-Verlag, Berlin, 2002.


\bibitem[MP10]{MP10}
P. Menal-Ferrer and  J. Porti, {\em Twisted cohomology for  hyperbolic three manifolds}, preprint (2010)

\bibitem[Mi66]{Mi66}
J. Milnor, {\em Whitehead torsion}, Bull. Amer. Math. Soc. 72 (1966), 358--426.


\bibitem[M\"u09]{Mu09}
W. M\"uller, {\em
Analytic torsion and cohomology of hyperbolic 3-manifolds}, Preprint of the Max-Planck Institute, Bonn (2009)

\bibitem[M\"u10]{Mu10}
W. M\"uller, {\em
The asymptotics of the Ray-Singer analytic torsion of hyperbolic 3-manifolds}, Preprint (2010)


\bibitem[Po97]{Po97}
J. Porti, {\em  Torsion de Reidemeister pour les vari\'et\'es hyperboliques},  Mem. Amer. Math. Soc. 128 (1997), no. 612

\bibitem[Ra10]{Ra10}
J. Raimbault, {\em Exponential growth of torsion in abelian coverings}, Preprint (2010)

\bibitem[Ri90]{Ri90}
R. Riley, {\em  Growth of order of homology of cyclic branched covers of knots}, Bull.
London Math. Soc. 22 (1990), 287-297.

\bibitem[Sc95]{Sc95}
K. Schmidt, {\em Dynamical Systems of Algebraic Origin}, Birkh\"auser Verlag, Basel, 1995.

\bibitem[Se10]{Se10}
M. H. \c{S}eng\"un, {\em On The Torsion Homology of Non-Arithmetic Hyperbolic Tetrahedral Groups}, Preprint (2010)


\bibitem[SW02]{SW02}
D. Silver and  S. Williams, {\em Mahler measure, links and homology
growth}, Topology 41 (2002), 979-991.

\bibitem[SW04]{SW04}
D. Silver and  S. Williams, {\em Mahler measure of Alexander polynomials}, J. London Math. Soc.
(2) 69 (2004), 767-782.
\bibitem[Tu01]{Tu01} V. Turaev, {\em Introduction to combinatorial torsions}, Birkh\"auser, Basel, (2001)

\end{thebibliography}
\end{document}